\expandafter\chardef\csname pre amssym.def at\endcsname=\the\catcode`\@ 
\catcode`\@=11 
 
\def\undefine#1{\let#1\undefined} 
\def\newsymbol#1#2#3#4#5{\let\next@\relax 
 \ifnum#2=\@ne\let\next@\msafam@\else 
 \ifnum#2=\tw@\let\next@\msbfam@\fi\fi 
 \mathchardef#1="#3\next@#4#5} 
\def\mathhexbox@#1#2#3{\relax 
 \ifmmode\mathpalette{}{\m@th\mathchar"#1#2#3}%
 \else\leavevmode\hbox{$\m@th\mathchar"#1#2#3$}\fi} 
\def\hexnumber@#1{\ifcase#1 0\or 1\or 2\or 3\or 4\or 5\or 6\or 7\or 8\or 
 9\or A\or B\or C\or D\or E\or F\fi} 
 
\font\tenmsa=msam10 
\font\sevenmsa=msam7 
\font\fivemsa=msam5 
\newfam\msafam 
\textfont\msafam=\tenmsa 
\scriptfont\msafam=\sevenmsa 
\scriptscriptfont\msafam=\fivemsa 
\edef\msafam@{\hexnumber@\msafam} 
\mathchardef\dabar@"0\msafam@39 
\def\dashrightarrow{\mathrel{\dabar@\dabar@\mathchar"0\msafam@4B}} 
\def\dashleftarrow{\mathrel{\mathchar"0\msafam@4C\dabar@\dabar@}} 
 
\def\ulcorner{\delimiter"4\msafam@70\msafam@70 } 
\def\urcorner{\delimiter"5\msafam@71\msafam@71 } 
\def\llcorner{\delimiter"4\msafam@78\msafam@78 } 
\def\lrcorner{\delimiter"5\msafam@79\msafam@79 } 
\def\yen{{\mathhexbox@\msafam@55 }} 
\def\checkmark{{\mathhexbox@\msafam@58 }} 
\def\circledR{{\mathhexbox@\msafam@72 }} 
\def\maltese{{\mathhexbox@\msafam@7A }} 
 
\font\tenmsb=msbm10 
\font\sevenmsb=msbm7 
\font\fivemsb=msbm5 
\newfam\msbfam 
\textfont\msbfam=\tenmsb 
\scriptfont\msbfam=\sevenmsb 
\scriptscriptfont\msbfam=\fivemsb 
\edef\msbfam@{\hexnumber@\msbfam}

\catcode`\@=\csname pre amssym.def at\endcsname 
 
\expandafter\ifx\csname pre amssym.tex at\endcsname\relax \else \endinput\fi 
\expandafter\chardef\csname pre amssym.tex at\endcsname=\the\catcode`\@ 
\catcode`\@=11 
\newsymbol\boxdot 1200 
\newsymbol\boxplus 1201 
\newsymbol\boxtimes 1202 
\newsymbol\square 1003 
\newsymbol\blacksquare 1004 
\newsymbol\centerdot 1205 
\newsymbol\lozenge 1006 
\newsymbol\blacklozenge 1007 
\newsymbol\circlearrowright 1308 
\newsymbol\circlearrowleft 1309 
\undefine\rightleftharpoons 
\newsymbol\rightleftharpoons 130A 
\newsymbol\leftrightharpoons 130B 
\newsymbol\boxminus 120C 
\newsymbol\Vdash 130D 
\newsymbol\Vvdash 130E 
\newsymbol\vDash 130F 
\newsymbol\twoheadrightarrow 1310 
\newsymbol\twoheadleftarrow 1311 
\newsymbol\leftleftarrows 1312 
\newsymbol\rightrightarrows 1313 
\newsymbol\upuparrows 1314 
\newsymbol\downdownarrows 1315 
\newsymbol\upharpoonright 1316 
  
\newsymbol\downharpoonright 1317 
\newsymbol\upharpoonleft 1318 
\newsymbol\downharpoonleft 1319 
\newsymbol\rightarrowtail 131A 
\newsymbol\leftarrowtail 131B 
\newsymbol\leftrightarrows 131C 
\newsymbol\rightleftarrows 131D 
\newsymbol\Lsh 131E 
\newsymbol\Rsh 131F 
\newsymbol\rightsquigarrow 1320 
\newsymbol\leftrightsquigarrow 1321 
\newsymbol\looparrowleft 1322 
\newsymbol\looparrowright 1323 
\newsymbol\circeq 1324 
\newsymbol\succsim 1325 
\newsymbol\gtrsim 1326 
\newsymbol\gtrapprox 1327 
\newsymbol\multimap 1328 
\newsymbol\therefore 1329 
\newsymbol\because 132A 
\newsymbol\doteqdot 132B 
  
\newsymbol\triangleq 132C 
\newsymbol\precsim 132D 
\newsymbol\lesssim 132E 
\newsymbol\lessapprox 132F 
\newsymbol\eqslantless 1330 
\newsymbol\eqslantgtr 1331 
\newsymbol\curlyeqprec 1332 
\newsymbol\curlyeqsucc 1333 
\newsymbol\preccurlyeq 1334 
\newsymbol\leqq 1335 
\newsymbol\leqslant 1336 
\newsymbol\lessgtr 1337 
\newsymbol\backprime 1038 
\newsymbol\risingdotseq 133A 
\newsymbol\fallingdotseq 133B 
\newsymbol\succcurlyeq 133C 
\newsymbol\geqq 133D 
\newsymbol\geqslant 133E 
\newsymbol\gtrless 133F 
\newsymbol\sqsubset 1340 
\newsymbol\sqsupset 1341 
\newsymbol\vartriangleright 1342 
\newsymbol\vartriangleleft 1343 
\newsymbol\trianglerighteq 1344 
\newsymbol\trianglelefteq 1345 
\newsymbol\bigstar 1046 
\newsymbol\between 1347 
\newsymbol\blacktriangledown 1048 
\newsymbol\blacktriangleright 1349 
\newsymbol\blacktriangleleft 134A 
\newsymbol\vartriangle 134D 
\newsymbol\blacktriangle 104E 
\newsymbol\triangledown 104F 
\newsymbol\eqcirc 1350 
\newsymbol\lesseqgtr 1351 
\newsymbol\gtreqless 1352 
\newsymbol\lesseqqgtr 1353 
\newsymbol\gtreqqless 1354 
\newsymbol\Rrightarrow 1356 
\newsymbol\Lleftarrow 1357 
\newsymbol\veebar 1259 
\newsymbol\barwedge 125A 
\newsymbol\doublebarwedge 125B 
\undefine\angle 
\newsymbol\angle 105C 
\newsymbol\measuredangle 105D 
\newsymbol\sphericalangle 105E 
\newsymbol\varpropto 135F 
\newsymbol\smallsmile 1360 
\newsymbol\smallfrown 1361 
\newsymbol\Subset 1362 
\newsymbol\Supset 1363 
\newsymbol\Cup 1264 
  
\newsymbol\Cap 1265 
  
\newsymbol\curlywedge 1266 
\newsymbol\curlyvee 1267 
\newsymbol\leftthreetimes 1268 
\newsymbol\rightthreetimes 1269 
\newsymbol\subseteqq 136A 
\newsymbol\supseteqq 136B 
\newsymbol\bumpeq 136C 
\newsymbol\Bumpeq 136D 
\newsymbol\lll 136E 
  
\newsymbol\ggg 136F 
  
\newsymbol\circledS 1073 
\newsymbol\pitchfork 1374 
\newsymbol\dotplus 1275 
\newsymbol\backsim 1376 
\newsymbol\backsimeq 1377 
\newsymbol\complement 107B 
\newsymbol\intercal 127C 
\newsymbol\circledcirc 127D 
\newsymbol\circledast 127E 
\newsymbol\circleddash 127F 
\newsymbol\lvertneqq 2300 
\newsymbol\gvertneqq 2301 
\newsymbol\nleq 2302 
\newsymbol\ngeq 2303 
\newsymbol\nless 2304 
\newsymbol\ngtr 2305 
\newsymbol\nprec 2306 
\newsymbol\nsucc 2307 
\newsymbol\lneqq 2308 
\newsymbol\gneqq 2309 
\newsymbol\nleqslant 230A 
\newsymbol\ngeqslant 230B 
\newsymbol\lneq 230C 
\newsymbol\gneq 230D 
\newsymbol\npreceq 230E 
\newsymbol\nsucceq 230F 
\newsymbol\precnsim 2310 
\newsymbol\succnsim 2311 
\newsymbol\lnsim 2312 
\newsymbol\gnsim 2313 
\newsymbol\nleqq 2314 
\newsymbol\ngeqq 2315 
\newsymbol\precneqq 2316 
\newsymbol\succneqq 2317 
\newsymbol\precnapprox 2318 
\newsymbol\succnapprox 2319 
\newsymbol\lnapprox 231A 
\newsymbol\gnapprox 231B 
\newsymbol\nsim 231C 
\newsymbol\ncong 231D 
\newsymbol\diagup 231E 
\newsymbol\diagdown 231F 
\newsymbol\varsubsetneq 2320 
\newsymbol\varsupsetneq 2321 
\newsymbol\nsubseteqq 2322 
\newsymbol\nsupseteqq 2323 
\newsymbol\subsetneqq 2324 
\newsymbol\supsetneqq 2325 
\newsymbol\varsubsetneqq 2326 
\newsymbol\varsupsetneqq 2327 
\newsymbol\subsetneq 2328 
\newsymbol\supsetneq 2329 
\newsymbol\nsubseteq 232A 
\newsymbol\nsupseteq 232B 
\newsymbol\nparallel 232C 
\newsymbol\nmid 232D 
\newsymbol\nshortmid 232E 
\newsymbol\nshortparallel 232F 
\newsymbol\nvdash 2330 
\newsymbol\nVdash 2331 
\newsymbol\nvDash 2332 
\newsymbol\nVDash 2333 
\newsymbol\ntrianglerighteq 2334 
\newsymbol\ntrianglelefteq 2335 
\newsymbol\ntriangleleft 2336 
\newsymbol\ntriangleright 2337 
\newsymbol\nleftarrow 2338 
\newsymbol\nrightarrow 2339 
\newsymbol\nLeftarrow 233A 
\newsymbol\nRightarrow 233B 
\newsymbol\nLeftrightarrow 233C 
\newsymbol\nleftrightarrow 233D 
\newsymbol\divideontimes 223E 
\newsymbol\varnothing 203F 
\newsymbol\nexists 2040 
\newsymbol\Finv 2060 
\newsymbol\Game 2061 
\newsymbol\mho 2066 
\newsymbol\eth 2067 
\newsymbol\eqsim 2368 
\newsymbol\beth 2069 
\newsymbol\gimel 206A 
\newsymbol\daleth 206B 
\newsymbol\lessdot 236C 
\newsymbol\gtrdot 236D 
\newsymbol\ltimes 226E 
\newsymbol\rtimes 226F 
\newsymbol\shortmid 2370 
\newsymbol\shortparallel 2371 
\newsymbol\smallsetminus 2272 
\newsymbol\thicksim 2373 
\newsymbol\thickapprox 2374 
\newsymbol\approxeq 2375 
\newsymbol\succapprox 2376 
\newsymbol\precapprox 2377 
\newsymbol\curvearrowleft 2378 
\newsymbol\curvearrowright 2379 
\newsymbol\digamma 207A 
\newsymbol\varkappa 207B 
\newsymbol\Bbbk 207C 
\newsymbol\hslash 207D 
\undefine\hbar 
\newsymbol\hbar 207E 
\newsymbol\backepsilon 237F 
\catcode`\@=\csname pre amssym.tex at\endcsname 
 
\magnification=1200 
\hsize=468truept 
\vsize=646truept 
\voffset=-10pt 
\parskip=1pc 
\baselineskip=14truept 
\count0=1 
 
\dimen100=\hsize 
 
\def\leftill#1#2#3#4{ 
\medskip 
\line{$ 
\vcenter{ 
\hsize = #1truept \hrule\hbox{\vrule\hbox to  \hsize{\hss \vbox{\vskip#2truept 
\hbox{{\copy100 \the\count105}: #3}\vskip2truept}\hss } 
\vrule}\hrule} 
\dimen110=\dimen100 
\advance\dimen110 by -36truept 
\advance\dimen110 by -#1truept 
\hss \vcenter{\hsize = \dimen110 
\medskip 
\noindent { #4\par\medskip}}$} 
\advance\count105 by 1 
} 
\def\rightill#1#2#3#4{ 
\medskip 
\line{ 
\dimen110=\dimen100 
\advance\dimen110 by -36truept 
\advance\dimen110 by -#1truept 
$\vcenter{\hsize = \dimen110 
\medskip 
\noindent { #4\par\medskip}} 
\hss \vcenter{ 
\hsize = #1truept \hrule\hbox{\vrule\hbox to  \hsize{\hss \vbox{\vskip#2truept 
\hbox{{\copy100 \the\count105}: #3}\vskip2truept}\hss } 
\vrule}\hrule} 
$} 
\advance\count105 by 1 
} 
\def\midill#1#2#3{\medskip 
\line{$\hss 
\vcenter{ 
\hsize = #1truept \hrule\hbox{\vrule\hbox to  \hsize{\hss \vbox{\vskip#2truept 
\hbox{{\copy100 \the\count105}: #3}\vskip2truept}\hss } 
\vrule}\hrule} 
\dimen110=\dimen100 
\advance\dimen110 by -36truept 
\advance\dimen110 by -#1truept 
\hss $} 
\advance\count105 by 1 
} 
\def\insectnum{\copy110\the\count120 
\advance\count120 by 1 
}

\font\ninerm=cmr9 
\font\eightrm=cmr8

\font\tenrm=cmr10 at 10pt 
 
\font\sc=cmcsc10

 
\def\msb{\fam\msbfam\tenmsb}

\def\bbc{{\msb C}}

\def\bbi{{\msb I}}

\def\bbp{{\msb P}} 
 
\def\bbr{{\msb R}}

\def\bbz{{\msb Z}}

\def\grD{\Delta}

\def\grL{\Lambda}

\def\grS{\Sigma}

\def\gro{\omega}

\def\la#1{\hbox to #1pc{\leftarrowfill}} 
\def\ra#1{\hbox to #1pc{\rightarrowfill}} 
 
\def\fract#1#2{\raise4pt\hbox{$ #1 \atop #2 $}}

\def\bowtie{\hbox to 1pt{\hss}\raise.66pt\hbox{$\scriptstyle{>}$} 
\kern-4.9pt\triangleleft} 
\def\hsmash{\triangleright\kern-4.4pt\raise.66pt\hbox{$\scriptstyle{<}$}} 
\def\boxit#1{\vbox{\hrule\hbox{\vrule\kern3pt 
\vbox{\kern3pt#1\kern3pt}\kern3pt\vrule}\hrule}}

\def\za{\vrule height6pt width4pt depth1pt}

\font\aa=eufm10

\def\Got#1{\hbox{\aa#1}}

\def\bfw{{\bf w}}

\def\cald{{\cal D}} 
 
\def\calf{{\cal F}}

\def\cals{{\cal S}}

\def\calz{{\cal Z}}

\def\gF{{\Got F}}

\def\Got#1{\hbox{\aa#1}}

\def\gsp1{{\Got s}{\Got p}(1)}

\font\svtnrm=cmr17

\font\bsc=cmcsc10 at 10truept

\def\endo{\hbox{End}}
\def\lcm{\hbox{lcm}}

\centerline{\svtnrm Rational Homology 5-Spheres}
\medskip
\centerline{\svtnrm with Positive Ricci Curvature}
\medskip

\bigskip
\centerline{\sc Charles P. Boyer~~ Krzysztof Galicki}
\footnote{}{\ninerm During the preparation of this work the authors 
were partially supported by NSF grant DMS-9970904.}
\bigskip
\bigskip

\centerline{\vbox{\hsize = 5.85truein
\baselineskip = 12.5truept
\eightrm
\noindent {\bsc Abstract:} We prove that for every integer $\scriptstyle{k>1}$ there is 
a simply connected rational homology
5-sphere $\scriptstyle{M^5_k}$ with spin such that
$\scriptstyle{H_2(M^5_k,\bbz)}$ has order $\scriptstyle{k^2},$ and $\scriptstyle{M^5_k}$ 
admits a Riemannian metric of positive Ricci curvature. Moreover, if 
the prime number decomposition of $\scriptstyle{k}$ has the form 
$\scriptstyle{k=p_1\cdots p_r}$ for distinct primes $\scriptstyle{p_i}$ then 
$\scriptstyle{M^5_k}$ is uniquely determined.}}
\tenrm
\bigskip

\bigskip
\baselineskip = 10 truept
\centerline{\bf  Introduction}  
\bigskip

Recently, Stephan Stolz brought to our attention the fact that preciously little is known 
about the existence of metrics of positive Ricci curvature on simply connected 
5-manifolds in the presence of torsion in $H_2.$ Indeed up until now there appears to be 
only one known non-trivial example of  a rational homology 5-sphere admitting metrics of 
positive Ricci curvature. This is somewhat surprising in light of the fact that it is precisely 
in dimension 5 that there is a diffeomorphism classification of compact simply connected 
5-manifolds.  In 1965 Barden proved the following remarkable theorem:

\noindent{\sc Theorem } [Bar]: \tensl The class of simply connected, closed, oriented, 
smooth, 5-manifolds is classifiable under diffeomorphism. Furthermore,
any such $M$ is diffeomorphic to one of the spaces
$M_{j;k_1,\ldots,k_s}=X_j\#M_{k_1}\#\cdots \#M_{k_s}$, where $-1\leq j\leq\infty,$
$s\geq0$, $1<k_1$ and $k_i$ divides $k_{i+1}$ or $k_{i+1}=\infty$. A complete set
of invariants is provided by $H_2(M,\bbz)$ and $i(M)$.
\tenrm

It is understood that $s=0$ means that no $M_{k_i}$ occurs. The diffeomorphism invariant 
$i(M)$ depends only on the 
second Stiefel-Whitney class $w^2(M)$. When $M$ is spin $i(M)=0$ as is
$w^2(M)=0$. Otherwise $i(M)\not=0$.
Barden's result is the extension of the similar theorem of Smale for
spin manifolds. In fact, in Barden's notation $M_{j;k_1,\ldots,k_s}$ is spin
if and only if $j=0$ in which case $X_0=S^5$ and
$M_{0;k_1,\ldots,k_s}=M_{k_1}\#\cdots \#M_{k_s}$ which is precisely
Smale's list [Sm]. Recall that here, 
$X_{-1}=SU(3)/SO(3)$ is the well-known symmetric space, 
$X_\infty$ is the non-trivial $S^3$ bundle over $S^2$.
For all other values of $j$, $X_j$ is a rational homology sphere 
determined by its 2-torsion which is $H_2(X_j,\bbz)=\bbz_{2^j}\oplus
\bbz_{2^j}$ and the value of $i(X_j)=j$. All $X_j$ other than $X_0$ are non-spin, as 
their second Stiefel-Whitney class $w^2(X_j)\not=0$. Now,
$M_k$, $1<k<\infty$ is a rational homology sphere with $H_2(M_k,\bbz)=\bbz_{k}\oplus
\bbz_{k}$  and vanishing $i(M_k)=0$. Finally, $M_\infty=S^2\times S^3$ is the trivial
$S^3$ bundle over $S^2$. All of the pieces in the above theorem are indecomposable 
with the exception of $X_1=X_{-1}\# X_{-1}$. Note that for $1<j<\infty,$ 
$M_{2^j}$ and $X_j$  have the same $H_2(M,\bbz)$
but $i(M_{2^j})=0$ while $i(X_j)=j$.

Regarding the question of positive Ricci curvature metrics on simply connected 
5-manifolds, there is the well-known result of
Sha and Yang [SY] which, in dimension 5, implies that any $k$-fold
connected sum of $S^2\times S^3$ admits a metric of positive Ricci
curvature (an alternative proof using methods similar to this note was given in [BGN2]).  
Using Barden's notation, there exists positive Ricci curvature
metrics on $M=M_\infty\#{\buildrel k\over \cdots}\# M_\infty$ for any
$k\geq1$. Now, $X_0\simeq S^5$ and the symmetric metric has constant curvature. 
Furthermore, $X_\infty$ admits a metric of positive Ricci curvature because of its
bundle structure [Na].
Remarkably, when the torsion in $H_2(M^5,\bbz)$ is non-trivial we know of just 
one example, namely $X_{-1}$ where  $H_2(X_{-1},\bbz)=\bbz_2$. However, our 
knowledge of this rational homology sphere is due to the special circumstance that it is the 
symmetric space $SU(3)/SO(3)$ which is known to admit a metric of positive Ricci 
curvature.  Indeed, the standard metric on this compact symmetric space 
happens to be a positive Einstein metric [Bes]. On the other
hand there are no known obstructions to the existence of positive Ricci curvature
metrics on any of the manifolds in Barden's classification. This poses
a natural 

\noindent{\sc Question}: \tensl Which smooth
simply connected closed 5-manifolds admit metrics of positive Ricci curvature?
\tenrm

We do not propose to answer this question here in its full generality.  Our
purpose is to prove the existence of metrics with positive Ricci 
curvature on an infinite number of the indecomposable building blocks in Barden's 
classification theorem by giving infinite families of simply 
connected rational homology 5-spheres $M^5_{k_i}$ admitting metrics with positive Ricci  
curvature. Our techniques use contact Riemannian geometry, more specifically Sasakian 
geometry.  Thus, our examples are complementary to the homogeneous manifold 
$X_{-1}=SU(3)/SO(3)$ whose second and third Stiefel-Whitney classes are both 
non-vanishing, and so it does not even admit an almost contact structure. In fact our 
method will not apply to any $M$ with $i(M)\not=0$. Our examples are necessarily spin 
and, in this note, they are rational homology spheres. By Smale's theorem [Sm] any rational 
homology 5-sphere $M^5$ with vanishing second Stiefel-Whitney class $w_2$ can be 
written uniquely as (for our purposes it is more convenient to rephrase  Smale's 
result in terms of elementary divisors instead of invariant factors as he does):
$$M^5=M^5_{p^{s_1}_1}\#\cdots \#M^5_{p^{s_r}_r}\leqno{(1)}$$
for some positive integers $r,s_1,\cdots,s_r$ where the $p_i$'s are (not necessarily 
distinct) primes, and  $H_2(M^5_{p^{s_i}_i},\bbz)=\bbz_{p^{s_i}_i}\oplus \bbz_{p^{s_i}_i}.$ 
A result of Geiges [Gei] says that any simply connected rational homology 5-sphere that 
admits a contact structure must be spin and thus a Smale manifold of the type given by 
equation (1). In [Gei] it was shown that a simply connected 5-manifold admits a contact 
structure if and only if its integral third Stiefel-Whitney class vanishes. Hence, from 
Barden's theorem, the $M_{j;k_1,\cdots,k_s}$ admits a contact structure if and only if 
$j=0$ or $1.$

Our main result is the following existence
theorem:

\noindent{\sc Theorem A}: \tensl For every integer $k>1$, there
exists a simply connected rational homology 5-sphere $M_k^5$ such that
$H_2(M_k^5,\bbz)$ has order $k^2$, $i(M_k^5)=0,$ and $M_k^5$ admits a Sasakian 
metric with positive Ricci curvature. 
\tenrm

Our construction gives all possible orders of $H_2,$ but it does not pin down the group 
precisely in all cases. However, the form of $H_2$ given by Smale's theorem says that the 
elementary divisors of $H_2$ must occur in pairs. Thus, when all of the primes $p_i$ in 
equation (1) are distinct and all of the $s_i$'s equal $1$, the order of $H_2$ uniquely 
determines the manifold. For example, if $|H_2(M^5,\bbz)|=36,$ the 
elementary divisors must be $\{2,2,3,3\}.$ This determines $M^5$ to be $M^5_2\#M^5_3.$
However, if for example, $|H_2(M^5,\bbz)|=64$ the elementary divisors can be 
$\{2^3,2^3\},\{2,2,2^2,2^2\},$ or $\{2,2,2,2,2,2\},$ giving the possibilities for $M^5$ as 
$M^5_{2^3},M^5_2\#M^5_{2^2},$ or $M^5_2\#M^5_2\#M^5_2.$ In this case we are 
unable to determine which manifolds occur. Thus, we stop short of proving that all 
simply connected rational homology 5-spheres with $i(M)=0$ admit metrics 
with positive Ricci curvature, although we certainly believe this to be the case. 

However, in cases when the order determines the group we have the following corollary

\noindent{\sc Corollary B}: \tensl For every positive integer $r$ and every list of distinct 
primes $p_1,\cdots, p_r,$ the manifolds
$$M^5=M^5_{p_1}\#\cdots \#M^5_{p_r}$$
admits Sasakian metrics with positive Ricci curvature. \tenrm

In the absence of known obstructions the question asked above is an intriguing one.
Since we are using methods of contact Riemannian geometry to
prove Theorem A we do not have any new ideas how to obtain such metrics
on manifolds with non-vanishing $i(M)$ such as for, for examples,
the rational homology 5-spheres $X_j$, $1\leq j<\infty$. Now,
$X_\infty$ is actually a contact, even a Sasakian manifold, but it is
not a positive Sasakian manifold as positivity implies spin [BGN2] and
$X_\infty$ is not spin. As already mentioned,
it follows from a theorem of Nash [Na] that $X_\infty$ 
does have a metric of positive Ricci curvature, but it cannot be Sasakian.

On the other hand, we believe that in the spin case one should be able to
extend Theorem A to many other cases. First it is reasonable that
actually all spin rational homology 5-spheres can be realized as links
with positive Sasakian structure. The main problem in proving this is
that we do not know how to compute the torsion group in all cases when it is
not determined by its order. One could also try to obtain examples of links $M$
when $H_2(M,\bbz)$ has both free part and torsion. Regarding this there is a conjectured 
algorithm due to Orlik [Or3] which has been verified in certain special cases, and we plan 
to address this question in the future. Another approach would be to use surgery, 
assuming that one has the existence of positive Ricci curvature metrics on the 
indecomposable pieces $M_k$. However, apparently the only known method, e.g. [SY], 
entails finding 2-spheres in $M_k$ such that the metric in a neighborhood of the $S^2$ is 
in standard form [St], and this appears to be obstructed by the 2-torsion.

\bigskip
\baselineskip = 10 truept
\centerline{\bf  1. Positive Sasakian Geometry}  
\bigskip

Let $(M,J)$ is a compact complex
manifold and $g$ a K\"ahler metric on $M$, with K\"ahler form M.
Suppose that $\rho'$ is a real, closed $(1,1)$-form on
$M$ with $[\rho']=2\pi c_1(M)$. Then there exists a unique K\"ahler
metric $g'$ on $M$ with K\"ahler form $\omega'$, such that
$[\omega]=[\omega']\in H^2(M,\bbr)$, and the Ricci form of $g'$ is
$\rho'$. The above statement is the celebrated Calabi Conjecture which
was posed by Eugene Calabi in 1954. The conjecture in its full generality was
eventually proved by Yau in 1976. In the Fano case when $c_1(M)>0$, i.e.,
when the first Chern class can be represented by a positive-definite
real, closed $(1,1)$-form $\rho'$ on $M$, the conjecture implies that
the K\"ahler form of $M$ can be represented by a metric of
positive Ricci curvature. The key idea behind the proof of Theorem A
is based on a more general Calabi Problem when $M$ is not necessarily
a smooth manifold but rather a $V$-manifold or an orbifold [DK, Joy]. 
In the contex of foliations one
actually can prove a ``transverse Yau theorem" and this was done
by El Kacimi-Alaoui in 1990 [ElK]. In [BGN2] we adapted this to a very special case of 
Sasakian foliations.

Recall [Bl, YK] that a Sasakian structure on a manifold $M$ of dimension $2n+1$
is a metric contact structure $(\xi,\eta,\Phi,g)$ such that the Reeb vector field $\xi$ is a
Killing field and whose underlying almost CR structure is integrable. Briefly,
let $(M,\cald)$ be a contact manifold, and choose a 1-form $\eta$ so that
$\eta\wedge (d\eta)^n\neq 0$ and $\cald=\ker~\eta.$  The pair $(\cald,\gro)$,
where $\gro$ is the restriction of $d\eta$ to $\cald$ gives $\cald$ the
structure of a symplectic vector bundle. Choose an almost complex structure
$J$ on $\cald$ that is compatible with $\gro,$ that is $J$ is a smooth section
of the endomorphism bundle $\endo~\cald$ that satisfies
$$ J^2= -\bbi, \qquad d\eta(JX,JY)=d\eta(X,Y), \qquad
d\eta(X,JX)>0\leqno{1.1}$$
for any smooth sections $X,Y$ of $\cald.$ Notice
that $J$ defines a Riemannian metric $g_\cald$ on $\cald$
by setting $g_\cald(X,Y) =d\eta(X,JY).$ One easily checks that $g_\cald$
satisfies the compatibility condition $g_\cald(JX,JY)=g_\cald(X,Y).$ Now we
can extend $J$ to an endomorphism $\Phi$ on all of $TM$ by putting $\Phi =J$
on $\cald$ and $\Phi\xi=0.$ Likewise we can extend the metric $g_\cald$ on
$\cald$ to a Riemannian metric $g$ on $M$ by setting
$$g=g_\cald + \eta\otimes \eta. \leqno{1.2}$$
The quadruple $(\xi,\eta,\Phi,g)$ is called a {\it metric contact structure}
on $M.$ If in addition $\xi$ is a Killing vector field and the almost complex
structure $J$ on $\cald$ is integrable the underlying almost contact
structure is said to be {\it normal} and $(\xi,\eta,\Phi,g)$ is called a
{\it Sasakian structure}. The fiduciary examples of compact Sasakian
manifolds are the odd dimensional spheres $S^{2n+1}$ with the standard contact
structure and standard round metric $g.$

Every Sasakian structure $\cals =(\xi,\eta,\Phi,g)$ has a 1-dimensional
foliation $\calf_\xi$ associated to it, defined by the flow of the Reeb vector
field $\xi$ and called the {\it characteristic foliation}. Associated with this foliation are 
important invariants, namely, the basic cohomology groups $H^p_B(\calf_\xi),$ (cf. [Ton]) 
and in particular one can consider the {\it basic first Chern class} $c_1(\calf_\xi)$ [ElK, 
BGN2] as an element in $H^2_B(\calf_\xi).$ These are not only invariants of the Sasakian 
structure, but of the entire deformation class of Sasakian structures. Notice that on a 
compact Sasakian manifold $H^2_B(\calf_\xi)\neq 0$ since $[d\eta]_B$ is a non-vanishing 
class. 

\noindent{\sc Definition} 1.3: \tensl A Sasakian manifold $M$ is said to
be {\it positive} if its basic first Chern class $c_1(\calf_\xi)$ can be
represented by a basic positive definite $(1,1)$-form. \tenrm

As in [BGN2, BGN4] we consider deformation classes $\gF(\calf_\xi)$ of Sasakian 
structures that have the same characteristic foliation. Recall that two Sasakian structures 
$\cals =(\xi,\eta,\Phi,g)$ and $\cals' =(\xi',\eta',\Phi',g')$ in $\gF(\calf_\xi)$ on 
a smooth manifold $M$ are said to be {\it $a$-homologous} if there is an $a\in 
\bbr^+$ such that $\xi'=a^{-1}\xi$ and $[d\eta']_B=a[d\eta]_B.$ On a rational homology 
sphere every $\cals\in \gF(\calf_\xi)$ belongs to precisely one of two $a$-homology 
classes corresponding to a given Sasakian structure or its conjugate.  In [BGN2] we 
proved, using El-Kacimi-Alaoui's [ElK] ``transverse Yau Theorem'',

\noindent{\sc Theorem} 1.4 [BGN2]: \tensl Let $\cals=(\xi,\eta,\Phi,g)$ 
be a positive Sasakian structure on a compact manifold $M$ of dimension 
$2n+1.$ Then $M$ admits a Sasakian structure $\cals'=(\xi',\eta',\Phi',g')$  
with positive Ricci curvature $a$-homologous to $\cals$ for some $a>0.$
\tenrm

Theorem 1.4 says that to prove the existence of a Sasakian metric with positive Ricci 
curvature it suffices to prove the existence of positive Sasakian structures.
In the next section we shall discuss how to construct positive Sasakian
structures on homotopy 5-spheres, and prove Theorem A of the Introduction.

\bigskip
\baselineskip = 10 truept
\centerline{\bf  2. The Construction}  
\bigskip

Our 5-manifolds are constructed as $k$-fold branched covers of $S^5$ branched over 
certain Seifert manifolds that are in turn $S^1$ orbifold  V-bundles over a compact 
Riemann surface of genus $g.$ Our construction is similar to that in [Sav]. Let 
$f_3(z_1,z_2,z_3)$ be a 
weighted homogeneous polynomial of an isolated hypersurface singularity in $\bbc^3$ with 
weights $\bfw=(w_1,w_2,w_3)$ and degree $d.$ The link $L_\bfw$ defined by 
$L_\bfw=\{f_3=0\}\cap S^5$ is a Seifert fibration over an algebraic curve $C_\bfw$ in the 
weighted projective space $\bbp(\bfw).$ Let $g=g(\bfw)$ denote the genus of the curve 
$C_\bfw.$ Then,

\noindent{\sc Proposition 2.1}: \tensl Let $L_f$ denote the link of the weighted 
homogeneous 
polynomial 
$$f=z_0^k+f_3(z_1,z_2,z_3)$$
with weights $(d,k\bfw)$ where $k$ is an integer $>1$ where $f_3$ is a weighted 
homogeneous polynomial of degree $d$ with weights $\bfw=(w_1,w_2,w_3)$ as above. 
Suppose further that $\gcd(d,k)=1.$ Then the link $L_f$ is a smooth simply connected 
rational  homology 5-sphere such that the order of 
$H_2(L_f,\bbz)$ is $k^{2g}.$ Furthermore, $L_f$  admits Sasakian metrics with positive 
Ricci curvature. \tenrm

\noindent{\sc Proof}: Let us briefly recall the construction of the Alexander 
(characteristic) polynomial $\grD_3(t)$ in [MO] associated to a 3-dimensional link 
$L_{f_3}.$ It is the characteristic polynomial of the monodromy map 
$\bbi-h_*:H_{2}(F,\bbz)\ra{1.3} H_{2}(F,\bbz)$ 
induced by the $S^1_\bfw$ action on the Milnor fibre $F.$ Thus, 
$\grD_3(t)=\det(t\bbi-h_*).$ Now both $F$ and its closure $\bar{F}$ are homotopy 
equivalent to a bouquet of $3$-spheres, and the boundary of $\bar{F}$ is the 
link $L_{f_3}.$ Now $L_{f_3}$ is connected and its Betti numbers  
$b_1(L_{f_3})=b_2(L_{f_3})$ equal the number of factors of $(t-1)$ in 
$\grD_3(t).$ Now since the curve $C_\bfw$ is algebraic, $b_1(L_{f_3})=2g$ where $g$ is 
the genus of $C_\bfw.$ 
Following Milnor and Orlik we let $\grL_j$ denote the divisor of $t^j-1$ in the group ring 
$\bbz[\bbc^*].$ Then the divisor of $\grD_3(t)$ is given by
$$\hbox{div}~\grD_3~= \prod_{i=1}^n({\grL_{u_i}\over v_i}-1)\leqno{(2.2)}$$
where we write ${d\over w_i}={u_i\over v_i}$ in irreducible form.
Using  the relations  $\grL_a\grL_b=\gcd(a,b)\grL_{lcm(a,b)},$ equation 
equation 1 takes the form
$\sum a_j\grL_j -1,$
where $a_j\in \bbz$ and the sum is taken over the set of all least common 
multiples of all combinations of the $u_1,\cdots,u_n.$ The Alexander 
polynomial is then given by
$$\grD_3(t) = (t-1)^{-1}\prod_j(t^j-1)^{a_j},\leqno{(2.3)}$$
and 
$$b_1(L_{f_3})~=~2g~=~\sum_j a_j-1. \leqno{(2.4)}$$

Now we compute the divisor $\hbox{div}~\grD$ of the Alexander polynomial $\grD_4(t)$ for 
$f.$ We have
$$\hbox{div}~\grD~=~(\grL_k-1)\hbox{div}~\grD_f~=~(\grL_k-1)\bigl(\sum 
a_j\grL_j-1\bigr)$$$$~=~\sum_j\gcd(k,j)a_j\grL_{\scriptstyle{\lcm(k,j)}} 
-\sum_ja_j\grL_j-\grL_k+1.$$
Since the $j$'s run through all the least common multiples of the set 
$\{u_1,\cdots,u_n\}$ and $\gcd(k,u_i)=1$ for all $i,$ we see that for all $j,$ 
$\gcd(k,j)=1.$ This implies 
$$b_2(L_f)~=~ \sum_ja_j-\sum_ja_j-1+1~=~0.$$
Thus, $L_f$ is a rational homology sphere. Next we compute the Alexander 
polynomial for $L_f.$
$$\grD_4(t) ~=~ {(t-1)\over (t^k-1)} \prod_j{(t^{kj}-1)^{a_j}\over   
(t^j-1)^{a_j}}$$$$=~(t^{k-1}+\cdots +t+1)^{-1} \prod_j \Bigl({t^{kj-1}+ 
\cdots +t+1\over t^{j-1}+\cdots +t+1}\Bigr)^{a_j}.\leqno{(2.5)}$$
This gives
$$\grD_4(1)~=~k^{-1}\prod_j\Bigl({kj\over j}\Bigr)^{a_j}~=~ k^{\grS_ja_j-1}~=~k^{2g}.$$
So by [MO] the order of $H_2(L_f,\bbz)$ is $\grD_4(1)=k^{2g}.$  

To finish the proof it suffices by Theorem 1.4 to show that the induced Sasakian structure 
on $L_f$ is positive, i.e.  that the basic first Chern class $c_1(\calf_\xi)\in 
H^2_B(\calf_\xi)$ is positive. Now from our previous work [BG1-2,BGN1-4] $L_f$ is the 
total 
space of a V-bundle over a K\"ahler orbifold $\calz_f.$ Moreover, $c_1(\calf_\xi)$ is just 
$c_1(\calz_f)$ pulled back to $L_f.$ Thus, it is enough to prove that $\calz_f$ is Fano. This 
follows from Lemma 2.6 below which is a special case of Lemma 3.12 of [BGN4].

\noindent{\sc Lemma} 2.6: \tensl As algebraic varieties $\calz_f$ is isomorphic to the 
weighted projective space $\bbp(\bfw).$ Hence, its Fano index is $|\bfw|=\sum_iw_i>0.$ 
\tenrm

\noindent This concludes the proof of Proposition 2.1. \hfill\za

Next we want to show that every order of $H_2$ can be realized. First we show that $f_3$ 
realizes curves of any genus. In fact there is a formula due to Orlik and Wagreich [OW] for 
the genus of the curve $C_\bfw$ which generalizes the well known genus formula for 
curves in $\bbp^2$ (See also the books [Dim, Or1]). It is
$$g(C_\bfw)={1\over 2}\Bigl({d^2\over w_1w_2w_3}-d\sum_{i<j}{\gcd(w_i,w_j)\over 
w_iw_j}+\sum_i{\gcd(d,w_i)\over w_i}-1\Bigr).\leqno{2.7}$$
We are interested in the case $g>0$  which implies $w_1+w_2+w_3\leq d$ [Or2].
It is easy to see that there are quasi-smooth weighted 
homogeneous polynomials $f_3$ with arbitrary genus, but we claim that genus one will 
suffice to realize all rational homology spheres of the form given in our main theorem.

\noindent{\sc Proposition 2.8}: \tensl For every integer $k>1$, there
exists a rational homology 5-sphere $M_k^5$ whose second homology group 
$H_2(M^5_k,\bbz)$ has order $k^2$ and that can be realized as the link $L_f$ of a 
weighted homogeneous polynomial $f$ given in Proposition 2.1 where $f_3$ cuts out a 
projective curve of genus one. Furthermore if $k$ has the form $k=p_1\cdots p_r$ for 
distinct primes $p_i,$ the manifold $M_k^5$ is uniquely determined up to diffeomorphism.  
\tenrm

\noindent{\sc Proof}: By Smale's classification theorem [Sm] and Proposition 2.1 it suffices 
to exhibit for each integer $k>1$ an infinite family of weighted homogeneous polynomials 
$f_3$ of prime degree $p$ with $g(C_\bfw)=1.$ For then for a given $k$ we can choose 
$p$ such that $\gcd(k,p)=1.$ The infinite family of polynomials is given by 
$$f_p(z_1,z_2,z_3)=z_1^p+z_2^2z_3+z_3^2z_1$$
with weights $\bfw=(1,{p+1\over 4},{p-1\over 2})$ and degree $p$ where $p$ is a prime of 
the form $p=4l-1.$ It is well known that there are an infinite number of such primes. The 
genus formula then gives
$$\eqalign{g=&{1\over 2}\Bigl({8p^2\over p^2-1} -p({4\over p+1} +{2\over p-1} +{8\over 
p^2-1}) +1+{4\over p+1} +{2\over p-1}-1\Bigr)\cr =&{1\over 2}{1\over 
p^2-1}(2p^2-6p+6p-2)=1,}$$
where we have used the fact that $\gcd({p+1\over 4},{p-1\over 2})=\gcd(l,2l-1)=1,$ and 
this proves the first statement. The second statement follows from the classification of 
finite Abelian groups and Smale's classification theorem [Sm].  \hfill\za
\bigskip
\noindent{\sc Remarks 2.7}: 
\item{(1)} In the case that $k=p_1\cdots p_r$ for distinct primes $p_i,$ the links cannot be 
realized using curves of higher genus $(g>1).$ 
\item{(2)} It is still an open question as to whether all simply connected rational homology 
5-spheres can be realized by our methods, and if so how does one distinguish the 
different elementary divisors. The curves of higher genus should play a role here.
\item{(3)} For each $k>1$ there are an infinite number of $p$'s that satisfy $\gcd(k,p)=1.$ 
This gives rise to an infinite number of Sasakian deformation classes $\gF(\calf_\xi)$ each 
with Sasakian metrics of positive Ricci curvature. It is also quite plausable that the 
different deformation classes belong to distinct underlying contact structures, but we have 
not proven this last statement.

\bigskip
\noindent{\sc Acknowledgments}: The authors would like to thank Stephan Stolz and 
Wolfgang Ziller
for fruitful discussions. The second author would also like to than
Max-Planck-Institut in Bonn for hospitality and support.

\bigskip
\centerline{\bf Bibliography}
\medskip
\medskip
\font\ninesl=cmsl9
\font\bsc=cmcsc10 at 10truept
\parskip=1.5truept
\baselineskip=11truept
\ninerm

\item{[Bar]} {\bsc D. Barden}, {\ninesl Simply-connected 5-manifolds}, Ann. of Math. 82 
(1965), 365-385.
\item{[Bes]} {\bsc A. Besse}, {\ninesl Einstein manifolds},
Springer-Verlag, Berlin-New York, 1987.
\item{[BG1]} {\bsc C. P. Boyer and  K. Galicki}, {\ninesl On Sasakian-Einstein
Geometry}, Int. J. Math. 11 (2000), 873-909.
\item{[BG2]} {\bsc C. P. Boyer and  K. Galicki}, {\ninesl New Einstein Metrics
in Dimension Five}, J. Diff. Geom. 57 (2001), 443-463.
\item{[BGN1]} {\bsc C. P. Boyer, K. Galicki, and M. Nakamaye}, {\ninesl On the
Geometry of Sasakian-Einstein 5-Manifolds}, submitted
for publication; math.DG/0012041.
\item{[BGN2]} {\bsc C. P. Boyer, K. Galicki, and M. Nakamaye}, {\ninesl On
Positive Sasakian Geometry}, submitted
for publication; math.DG/0104126.
\item{[BGN3]} {\bsc C. P. Boyer, K. Galicki, and M. Nakamaye}, {\ninesl Einstein Metrics 
on Rational Homology 7-Spheres}, submitted for publication; math.DG/0108113.
\item{[BGN4]} {\bsc C. P. Boyer, K. Galicki, and M. Nakamaye}, {\ninesl Sasakian Geometry, 
Homotopy Spheres and Positive Ricci Curvature}, submitted for publication.
\item{[Bl]} {\bsc D.E. Blair}, {\ninesl Contact Manifolds in Riemannian Geometry}, LNM 509,
Springer-Verlag, 1976.
\item{[DK]} {\bsc J.-P. Demailly and J. Koll\'ar}, {\ninesl Semi-continuity of
complex singularity exponents and K\"ahler-Einstein metrics on Fano
orbifolds}, preprint AG/9910118, Ann. Scient. Ec. Norm. Sup. Paris 34 (2001), 525-556.
\item{[Dim]} {\bsc A. Dimca}, {\ninesl Singularities and Topology of
Hypersurfaces}, Springer-Verlag, New York, 1992.
\item{[ElK]} {\bsc A. El Kacimi-Alaoui}, {\ninesl Op\'erateurs transversalement
elliptiques sur un feuilletage riemannien et applications}, Compositio
Mathematica 79 (1990), 57-106.
\item{[Gei]} {\bsc H. Geiges}, {\ninesl Contact Structures on 1-connected 5-manifolds}, 
Mathematika 38 (1991), 303-311.
\item{[Joy]} {\bsc  D. Joyce}, {\ninesl
Compact manifolds with special holonomy},
Oxford Mathematical Monographs, Oxford University Press, Oxford 2000.
\item{[MO]} {\bsc J. Milnor and P. Orlik}, {\ninesl Isolated singularities
defined by weighted homogeneous polynomials}, Topology 9 (1970), 385-393.
\item{[Na]} {\bsc J. Nash}, {\ninesl Positive Ricci curvature on fibre bundles}, J. Diff. 
Geom. 14 (1979), 241-254.
\item{[Or1]} {\bsc P. Orlik}, {\ninesl 
Seifert manifolds},
Lecture Notes in Mathematics, Vol. 291, Springer-Verlag, Berlin-New York, 1972.
\item{[Or2]} {\bsc P. Orlik}, {\ninesl Weighted homogeneous polynomials and
fundamental groups}, Topology 9 (1970), 267-273. 
\item{[Or3]} {\bsc P. Orlik}, {\ninesl On the homology of weighted homogeneous
manifolds}, Proc. 2nd Conf. Transformations Groups I, LNM 298,
Springer-Verlag, (1972), 260-269.
\item{[OW]} {\bsc P. Orlik and  P. Wagreich}, {\ninesl 
Isolated singularities of algebraic surfaces with $\bbc^*$ action},
Ann. of Math. (2) 93 (1971) 205--228.
\item{[Sav]} {\bsc I.V. Savel'ev} {\ninesl Structure of Singularities of a 
Class of Complex Hypersurfaces}, Mat. Zam. 25 (4) (1979) 497-503; English 
translation: Math. Notes 25 (1979), no. 3--4, 258--261.
\item{[Sm]} {\bsc S. Smale}, {\ninesl On the structure of 5-manifolds},
Ann. Math. 75 (1962), 38-46.
\item{[St]} {\bsc S. Stolz}, private communication.
\item{[SY]} {\bsc J.-P. Sha and D.-G Yang}, {\ninesl Positive Ricci
curvature on the connected sums of $S^n\times S^m$}, J. Diff. Geom. 33
(1991), 127-137.
\item{[Ton]} {\bsc Ph. Tondeur}, {\ninesl Geometry of Foliations}, Monographs
in Mathematics, Birkh\"auser, Boston, 1997.
\item{[YK]} {\bsc K. Yano and M. Kon}, {\ninesl
Structures on manifolds}, Series in Pure Mathematics 3,
World Scientific Pub. Co., Singapore, 1984.
\medskip
\bigskip \line{ Department of Mathematics and Statistics
\hfil February 2002} \line{ University of New Mexico \hfil }
\line{ Albuquerque, NM 87131 \hfil } \line{ email: cboyer@math.unm.edu,
galicki@math.unm.edu \hfill} \line{ web pages:
http://www.math.unm.edu/$\tilde{\phantom{o}}$cboyer,
http://www.math.unm.edu/$\tilde{\phantom{o}}$galicki \hfil}

\bye